\documentclass[11pt,a4paper]{amsart}
\usepackage{amssymb}
\usepackage{amsmath}
\usepackage{amsfonts,a4wide}
\usepackage{bbm}
\usepackage{caption}
\usepackage{enumerate}
\usepackage[latin1]{inputenc}
\usepackage{shadow}
\usepackage[all]{xy}

\newtheorem{thm}[subsection]{Theorem}

\newtheorem{prop}[subsection]{Proposition}

\newtheorem{lemma}[subsection]{Lemma}
\newtheorem{conj}[subsection]{Conjecture}
\newtheorem{remark}[subsection]{Remark}
\newcommand{\cat}{\mathcal }
\newcommand{\A}{\cat A}
\newcommand{\llra}[1]{\stackrel{#1}{\longrightarrow}}
\newcommand{\HP}{\mathbb{H}{\rm P}}
\newcommand{\HPi}{\HP^\infty}
\newcommand{\HPbar}{\overline{\HPi}}
\newcommand{\tens}{\otimes}
\newcommand{\Z}{\mathbb{Z}}
\newcommand{\R}{\mathbb{R}}
\newcommand{\Q}{\mathbb{Q}}
\newcommand{\F}{\mathbb{F}}
\newcommand{\Fbar}{\overline{F_3}}
\newcommand{\ra}{\rightarrow}
\newcommand{\lla}{\longleftarrow}
\newcommand{\lra}{\longrightarrow}

\DeclareMathOperator{\Ext}{Ext}
\DeclareMathOperator{\Hom}{Hom}
\DeclareMathOperator{\hofib}{hofib}

\begin{document}

\title{Self maps of $\HP^n$ via the unstable Adams spectral sequence}
\author{Gustavo Granja}
\address{ Departamento de Matem\'atica, 
Instituto Superior T\'ecnico, Universidade de Lisboa,
Av. Rovisco Pais, 1049-001 Lisboa, Portugal}
\email{ggranja@math.tecnico.ulisboa.pt}
\date{\today}

{\abstract We use obstruction theory based on the unstable Adams
spectral sequence to construct self maps of finite quaternionic projective 
spaces. As a result, a conjecture of Feder and Gitler regarding the 
classification of self maps up to homology is proved in two new
cases.}

\maketitle

\section{Introduction}

Let $\HP^n$ denote $n$-dimensional quaternionic projective space and recall
that $H^*(\HP^n;\Z)=\Z[u]/u^{n+1}$ with $|u|=4$. Given a self 
map $f$ of $\HP^n$, the \emph{degree} of $f$ is the integer 
$\deg(f)$ such that $f^*(u)=\deg(f)u$.

The homotopy classification of self maps of $\HP^\infty$ is well known:
 self maps are classified by their degrees and the allowable
degrees are zero and the odd square integers \cite{Ms}. The situation for
 finite projective spaces is more complicated. It is not true in general
that self maps are classified by their degrees (see \cite{MR,IMO,GS}), and even the 
set of possible degrees is unknown. We will write 
\[\mathbf{R}_n=\{ \deg(f) \in \Z \colon f: \HP^n \to \HP^n\}\] 
for the set of possible degrees. Alternatively, 
$\mathbf{R}_n$ can be described as the set of homotopy classes of self maps of $S^3$ which 
admit an $A_n$ structure (see \cite{St}).

There is a conjectural description of $\mathbf{R}_n$ due to Feder and Gitler. 
Consider the congruences
\begin{equation}
\label{cong}
C_n: \hspace{1cm} \prod_{i=0}^{n-1} (k-i^2) \  \equiv 0 \bmod \left\{
\begin{array}{ll}
(2n)! & \text{ if } n \text{ is even,} \\
 & \\
\frac{(2n)!}{2} & \text{ if } n \text{ is odd.} 
\end{array}
\right.
\end{equation}
and define the sets of integers
\[\mathbf{FG}_n=\{k \in \Z \quad | \quad k \text{ is a solution of } C_1,\ldots,C_n\}.\]
In \cite{FG}, Feder and Gitler computed the restrictions imposed by $K$-theory and Adams operations 
on the degrees of self maps of $\HP^n$: they are precisely that the degree should be in $\mathbf{FG}_n$.
Therefore $\mathbf{R}_n \subset \mathbf{FG}_n$. 
\begin{conj}[Feder and Gitler]
\label{conje}
$\mathbf{R}_n = \mathbf{FG}_n$.
\end{conj}
This is trivial for $n=1$. A proof for $n=2$ is
contained in \cite{AC}, and for $n=3$ in \cite{MG2}. Feder and Gitler showed in
\cite{FG} that the integers in $\mathbf{FG}_\infty$ are precisely the odd square integers
and zero so the conjecture also holds for $n=\infty$. In fact, self maps of $\HP^\infty$ with
 these degrees were first constructed by Sullivan in \cite{Su}, motivating 
Feder and Gitler's work.
McGibbon has also proved that the conjecture holds stably, in a suitable sense (see
\cite[Theorem 3.5]{MG1}).

The solutions to the congruences \eqref{cong} are more easily described
one prime at a time. It turns out (see Proposition \ref{localcong}) that 
the $p$-local integers which satisfy $C_1, \ldots, C_n$ consist of 
the $p$-adic squares (if $p=2$ they must also be units or zero) together
with all multiples of $p^{e(p,n)}$ for a certain function $e$ described in 
\eqref{exponent}. The former can all be realized as degrees of self maps of the 
$p$-localization $\HPi_{(p)}$ by work of Rector (see Proposition \ref{rector}), so 
to prove the conjecture it is enough to build self maps of 
$\HP^n_{(p)}$ with degree an arbitrary multiple of $p^{e(p,n)}$. 

The obstruction to extending a self map of $\HP^n$ to $\HP^{n+1}$
is an element in $\pi_{4n+2}S^3$ (see \eqref{obst}). 
We will use some of the information available on $\pi_*S^3$ to  
analyse this obstruction and produce self maps of $\HP^n_{(p)}$ with degree any multiple
 of $p^{f(p,n)}$ where  $f$ is a different function described in \eqref{defnf} (see Theorem \ref{main2}). 
In general, $e(p,n)<f(p,n)$ but for $n\leq 5$ the two functions agree which implies the following result.
\begin{thm}
\label{mainthm}
The Feder-Gitler conjecture holds for $n\leq 5$.
\end{thm}
This is qualitatively stronger evidence for the conjecture than had previously been obtained as 
the obstruction group $\pi_{4n-2}S^3$ is detected by $K$-theory and Adams operations
 (i.e. by the Adams  $e$-invariant) only for $n\leq 3$. For $n\geq 3$ it is easy to produce 
self maps of $\HP^n$ with degree in $\mathbf{FG}_{n+1}$ which do not extend one stage. 
Theorem \ref{mainthm} is proved by inductively
building maps of high Adams filtration thus ensuring their obstructions to extension 
are detected by the $e$-invariant. 

The paper is organized as follows. In section \ref{sectionlocal}  we describe 
the $p$-local solutions to the congruences \eqref{cong} and quote results of 
Rector which reduce the Feder-Gitler conjecture to producing maps with 
degree any multiple of $p^{e(p,n)}$.
In section \ref{sectionadams}, we gather the necessary information on the unstable Adams
spectral sequence for $\HP^\infty$.
In section \ref{sectionmaps} we show that the degree of a self map of $\HP^n$
is in $\mathbf{FG}_{n+1}$ if and only if the Adams $e$-invariant 
of its obstruction to extension vanishes. The results quoted in section \ref{sectionadams} then allow
us to construct self maps of $\HP^n_{(p)}$ of degree any multiple of $p^{f(p,n)}$ and
thus  prove Theorem \ref{mainthm}.

\subsection{Acknowledgments}  The results of this paper are part of the author's 1999 MIT Ph.D. Thesis
written under the supervision of Mike Hopkins whom the author thanks for his guidance. This research
was partially supported by grants from FCT-Portugal, FLAD and Fundação Gulbenkian. The author
warmly thanks John R. Harper for pointing out a mistake in a previous version of this paper and 
explaining how to fix it.

\section{The local version of the conjecture}
\label{sectionlocal}

We will write $|k|_p=\sup\{l \colon p^l | k \}$ for the 
$p$-adic valuation of $k$ and $\lfloor x \rfloor$ for the greatest integer $\leq x$.
In this section we reduce the Feder-Gitler conjecture to constructing
self maps of the localizations $\HP^n_{(p)}$ of degree $k$ for all $k$ with
$|k|_p \geq e(p,n)$ for a certain function $e(p,n)$.

\smallskip

First note that, by cellular approximation, the inclusion 
of $\HP^n$ in $\HP^\infty$ induces a 
bijection $[\HP^n, \HP^n] \to [\HP^n,\HP^\infty]$. As 
$\HP^\infty_{(0)} = K(\Q,4)$, it follows easily from 
Sullivan's arithmetic square that in order to produce a 
self map of $\HP^n$ of degree $k$, it suffices to do so after 
localizing at each prime $p$ (or, in fact, at each prime $p\leq 2n-1$).

We will write $\mathbf{R}_{n,p}=\{\deg(f) \in \Z_{(p)} | 
f:\HP^n_{(p)} \to \HP^n_{(p)}\}$ and $\mathbf{FG}_{n,p}$ for the set 
of solutions in $\mathbf{\Z}_{(p)}$ of the congruences \eqref{cong}. 
Let $\mathbf{D}_p$ denote the set defined by

\[ \mathbf{D}_p = \left\{ \begin{array}{ll}
\{ k \in \Z_{(2)} \ | \ k=0 \text { or } k \equiv 1 \bmod 8 \} &
 \text{ if } p=2, \\
\{k \in \Z_{(p)} \ | \ k = u^2 \text{ for some } u \in \hat{\Z}_p \}
& \text{ if } p>2.
\end{array}
\right.
\]
It will soon emerge that $\mathbf{D}_p=\mathbf{R}_{\infty,p}=
\mathbf{FG}_{\infty,p}$. Note that $\mathbf{D}_2$ consists of  
the $2$-adic unit squares in $\Z_{(2)}$ together with $0$. We define the function
\begin{equation}
\label{exponent}
e(p,n) =  \begin{cases}
0 & \text{ if } n=1, \\
\#\{k\leq n \mid k=p^j \text{ or } 
k=(\frac{p+1}{2})p^{j-1} \text{ for some } j\geq 1 \}
& \text{ if } p>2,\\
1+2\lfloor \log_2(n)\rfloor & \text{otherwise.} \\
\end{cases}
\end{equation}

\begin{prop}
\label{localcong}
Let $p$ be a prime and $n<\infty$. Then
\[ \mathbf{FG}_{n,p}= \mathbf{D}_p \cup p^{e(p,n)}\Z_{(p)}.\]
\end{prop}
\begin{proof}
The sets $\mathbf{D}_{p}$ are the closure of 
$\mathbf{R}_\infty=\{0,1,9,25,\ldots\}$ in $\Z_{(p)}$ with 
the $p$-adic topology. Since the solution set of each 
congruence is closed we conclude that $\mathbf{D}_p \subset
\mathbf{FG}_{n,p}$ for every $n,p$. We will conclude the proof
in the case $p=2$ (the case when $p$ is odd is similar).

If $k\in \Z_{(2)}$ is a solution of $C_2$ then either $k\equiv 1 \bmod 8$,  
or $k=2^3l$ for some $l\in \Z_{(2)}$ so it suffices to consider 
solutions $k$ with $|k|_2\geq 3$. For these, since $|(4m+2)!/2|_2=|(4m)!|_2$ it 
suffices to consider the congruences $C_{2m}$. So let $|k|_2\geq 3$.
Factoring out units in $\mathbb{Z}_{(2)}$ we see that $k$ is a solution of $C_{2m}$
if and only if  
\[
\prod^{m-1}_{i=0} (k-(2i)^{2}) \equiv 0 \bmod 2^{|(4m)!|_2}.
\]
Writing $k=4l$ and noting that $|(4m)!|_2=2m+|(2m)!|_2$, this becomes
\[
\prod^{m-1}_{i=0} (l-i^{2}) \equiv 0 \bmod 2^{|(2m)!|_2} 
\]
so if $k=4l \in \mathbf{FG}_{2m,2}$ then $l \in \mathbf{FG}_{m,2}$. 
Since $l \not\in \Z_{(2)}^{\times}$, assuming the 
result is true for all $n<2^d$ and letting $2^d\leq 2m<2^{d+1}$, it follows by induction
that $\mathbf{FG}_{2m,2}-\mathbf{D}_2 \subset 2^{2d+1}\Z_{(2)}$. 

It remains to check that for $2^d \leq 2m < 2^{d+1}$ we have $k \in 
\mathbf{FG}_{2m,2}$ as long as $|k|_2\geq 2d+1$. That is, for $j<2^d$ and $|k|_2 \geq 2d+1$,
the following equivalent inequalities hold 
\begin{eqnarray*}
  \left| \prod_{i=0}^{2j-1} (k-i^2) \right|_2 &\geq &|(4j)!|_2 \\
 |k|_2 +\sum_{i=1}^{2j-1} 2|i|_2 & \geq & 2j + j + \ldots + \lfloor j/2^{d-1} \rfloor \\
 |k|_2 + 2(j + \ldots + \lfloor j/2^{d-1}\rfloor ) - 2|2j|_2  & \geq & 2j + j + \ldots + \lfloor j/2^{d-1} \rfloor\\
 |k|_2 + j + \ldots + \lfloor j/2^{d-1}\rfloor  & \geq & 2j + 2|2j|_2.
\end{eqnarray*}
The last inequality is easy to check writing $2j=2^rv$ with $v$ odd, noticing that $r\leq d$ and $v<2^{d-r+1}$.
\end{proof}
As $e(p,n)$ tends to $\infty$ with $n$, this means that $\mathbf{FG}_{\infty,p} = \mathbf{D}_p$.
The following result of Rector implies that in fact we have 
$\mathbf{R}_{\infty,p}=\mathbf{FG}_{\infty,p} = \mathbf{D}_p$ and therefore, that in order to prove 
Conjecture \ref{conje} it suffices to construct self maps of $\HP^n$ with degrees all multiples of $p^{e(p,n)}$.
We note that, although the statement on \cite[p.103]{Re} in the case when $p$ is odd 
 concerns only the units in $\mathbf{D}_p$, the argument given applies verbatim to all  of $\mathbf{D}_p$.
\begin{prop}[Rector]
\label{rector}
If $k \in \mathbf{D}_p$ then there is a self map of $\HP^\infty_{(p)}$ of degree $k$.
\end{prop}

\section{The unstable Adams spectral sequence for $\HP^\infty$}
\label{sectionadams}

Let $\cat U$ denote the category of unstable modules over the mod $p$ Steenrod
algebra. If $M \in \cat U$, $U(M)$ denotes the free unstable
algebra generated by $M$ and for  $P \in \cat U$ a free module, we write 
$K(P)$ for the product of mod $p$ Eilenberg-Maclane spaces such that $H^*(K(P);\F_p) = U(P)$.
Recall from \cite{HM2} (see also \cite[Section 7.2]{Ha})  that if $X$ is a simply connected space such that
$H^*(X;\F_p) = U(M)$ for some $M \in \cat U$ of finite type, then given a free resolution 
\[ 0 \lla M \lla P_0 \lla P_1 \lla P_2 \lla \ldots \]
of $M$ with $P_s$ of finite type and $0$ in degree $\leq s+1$, one can construct a tower
 of principal fibrations over $X$, called an \emph{Adams resolution for $X$}
\begin{equation}
\label{tower}
\xymatrix{
  {\vdots} \ar[d]^{p_2} &\\
  X\{2\} \ar[r]^{j_2} \ar[d]^{p_1} & K(\Omega^2 P_2) \\
  X\{1\} \ar[r]^{j_1} \ar[d]^{p_0} & K(\Omega P_1) \\
  X\{0\}=X \ar[r]^{j_0} & K(P_0)
}
\end{equation}
where $\Omega$ denotes the algebraic loop functor (the left adjoint to suspension).
The corresponding homotopy spectral sequence is the Massey-Peterson version of the unstable
Adams spectral sequence. For a simply connected space of finite type this spectral 
sequence converges strongly and
\[ E_2^{s,t}(X) = \Ext_{\mathcal{U}}^s(M,\Sigma^t\mathbb{F}_p) \Rightarrow  
\pi_{t-s}X\tens\widehat{\Z}_p. \]
The \emph{Adams filtration} of a map $f:W \ra X$ is the largest $s$
such that $f$ factors through $X\{s\}$. It's easy to check that $f$ has Adams filtration $s$ if and only if it can be
factored as a composite of $s$ maps (of spaces) inducing $0$ on mod $p$ homology, but not as a 
composite of $s+1$ such maps. The following result is an immediate consequence of the 
definition of the homotopy spectral sequence.
\begin{lemma}
\label{coninj}
Let $X$ be a simply connected $p$-local space of finite type with $H^*(X;\F_p)=U(M)$ and
let $E_r^{s,t}(X) \Rightarrow \pi_{t-s} X \otimes \widehat{\Z}_p$ denote the Massey-Peterson
spectral sequence associated to a minimal resolution of $M$ (so that $d_1=0$). Then
\begin{enumerate}[(a)]
\item  If $E_2^{s,n+s}=0$ for $s\geq r$, then $\pi_n X\{r\}=0$.
\item  Suppose $0\leq r < s$ and the differentials $d_k: E_k^{u,n+1+u} \to E_k^{u+k,n+u+k}$ are $0$
for $r\leq u<s$ and $k\geq 2$. Then the map $\pi_n X\{s\} \to \pi_n X\{r\}$ is injective.
\end{enumerate}
\end{lemma}

The algebra $H^*(\HPi;\F_p)$ is only of the form $U(M)$ when $p=2$, but as pointed out to me 
by John R. Harper, this problem can be overcome by replacing $\HPi$ with the homotopy fiber 
$$ \HPbar = \hofib \left( \HPi \xrightarrow{\alpha} K(\Z/p,4) \right) $$
where $\alpha$ denotes a generator of $H^4(\HPi;\F_p)$ (when $p=2$ and we use a minimal resolution,
$\HPbar = \HPi\{1\}$).
Indeed, it follows easily from the Serre spectral sequence of the fibration
$$ K(\Z/p,3) \to \HPbar \to \HPi $$
that, denoting by $F_3$ the free unstable module on a class of degree $3$ and by $\Fbar$
 the kernel of the canonical map $F_3 \to \Sigma^3\F_p$, we have  $H^*(\HPbar;\F_p) = U(\Fbar)$.
It follows from the long exact sequence in $\Ext$ that the the $E_2$ terms of the 
Massey-Peterson spectral sequences for $S^3$ and $\HPbar$ are related
by 

\begin{equation}
\label{relE2}
E_2^{s,t}(\HPbar) = \Ext^s_{\cat U}( \Fbar, \Sigma^t \F_p) \cong \Ext^{s+1}_{\cat U}(\Sigma^3 \F_p , \Sigma^t F_p) =
E_2^{s+1,t}(S^3).
\end{equation}

We will now invoke some well known information about the spectral sequence for $S^3$ (cf. \cite{Ta, CM}),
from which it will follow that Figures \ref{v1even} and \ref{v1odd} describe the portion of the $E_2$ term 
of the spectral sequence for $\HPbar$ along the vanishing line. 
Vertical lines represent multiplication by $p$ and the slanted lines
composition with $\eta$ if the prime is $2$ and with
 $\alpha_1$ if the prime is odd. If $p$ is an odd prime, we set $q=2(p-1)$.

\begin{figure}[h]
\begin{center}
\newcounter{c}
\setlength{\unitlength}{.7mm}
\begin{picture}(200,55)(0,0)
\put(2,40){\makebox{$s$}}
\put(10,10){\vector(0,1){30}}
\setcounter{c}{0}
\multiput(3,8.5)(0,5){6}{
   \makebox(0,0)[b]{\arabic{c}}
   \addtocounter{c}{1}
}
\multiput(10,10)(0,5){6}{\line(1,0){1}}
\put(10,10){\vector(1,0){50}}
\multiput(15,10)(5,0){8}{\line(0,1){1}}
\put(65,5){\makebox{$t-s$}}
\setcounter{c}{4}
\multiput(14,5)(5,0){7}{
   \makebox(0,0)[b]{\arabic{c}}
   \addtocounter{c}{1}
}
\put(15,10){\line(0,1){25}}
\multiput(15,10)(0,5){6}{\circle{1}}
\put(20,10){\line(1,1){10}}
\put(18,12){$\eta$}
\multiput(20,10)(5,5){3}{\circle{1}}
\put(30,15){\circle{1}}
\put(30,15){\line(0,1){5}}
\put(30,15){\line(1,1){10}}
\put(35,20){\circle*{1}}
\put(40,25){\circle*{1}}
\multiput(42,12)(3,0){4}{\line(1,0){2}}

\put(85,43){\makebox{$s=4k+$}}
\put(110,10){\vector(0,1){35}}
\setcounter{c}{-3}
\multiput(103,12)(0,5){6}{
   \makebox(0,0)[b]{\arabic{c}}
   \addtocounter{c}{1}
}
\multiput(110,13.5)(0,5){6}{\line(1,0){1}}
\put(93,0){\makebox{$t-s=8k+$}}
\put(120,5){\vector(1,0){50}}
\setcounter{c}{4}
\multiput(122,0)(5,0){8}{
   \makebox(0,0)[b]{\arabic{c}}
   \addtocounter{c}{1}
}
\multiput(123,5)(5,0){8}{\line(0,1){1}}
\put(123,18.5){\circle{1}}
\put(128,23.5){\circle{1}}
\put(133,28.5){\circle*{1}}
\put(133,23.5){\circle*{1}}
\put(128,28.5){\circle{1}}
\put(133,33.5){\circle{1}}
\put(138,38.5){\circle{1}}
\put(138,33.5){\circle{1}}
\put(138,28.5){\circle*{1}}
\put(143,38.5){\circle*{1}}
\put(143,33.5){\circle*{1}}
\put(148,43.5){\circle*{1}}
\put(123,18.5){\line(1,1){10}}
\put(128,28.5){\line(1,1){10}}
\put(133,23.5){\line(1,1){10}}
\put(138,33.5){\line(1,1){10}}
\put(133,23.5){\line(0,1){5}}
\put(138,33.5){\line(0,1){5}}
\multiput(148,28.5)(3,0){4}{\line(1,0){2}}
\end{picture}
\caption{$E_2^{(s,t)}(\HPbar)$ for $p=2$.} \label{v1even}
\end{center}

\begin{center}
\setlength{\unitlength}{.7mm}
\begin{picture}(200,45)(0,0)
\put(5,30){\makebox{$s$}}
\put(10,10){\vector(0,1){25}}
\setcounter{c}{0}
\multiput(3,8.5)(0,5){4}{
   \makebox(0,0)[b]{\arabic{c}}
   \addtocounter{c}{1}
}
\multiput(10,10)(0,5){4}{\line(1,0){1}}
\put(10,10){\vector(1,0){50}}
\multiput(15,10)(5,0){8}{\line(0,1){1}}
\put(65,5){\makebox{$t-s$}}
\put(14,5){\makebox{4}}
\put(15,10){\line(0,1){15}}
\multiput(15,10)(0,5){4}{\circle{1}}
\put(30,10){\circle{1}}
\put(27,12){$\alpha_1$}
\put(45,15){\circle*{1}}
\put(50,15){\circle{1}}
\put(30,10){\line(3,1){15}}
\put(25,5){\makebox{$3+q$}}
\put(40,5){\makebox{$2+2q$}}

\put(88,30){\makebox{$s=k+$}}
\put(110,10){\vector(0,1){25}}
\setcounter{c}{-2}
\multiput(103,12)(0,5){3}{
   \makebox(0,0)[b]{\arabic{c}}
   \addtocounter{c}{1}
}
\multiput(110,13.5)(0,5){3}{\line(1,0){1}}
\put(92,0){\makebox{$t-s=qk+$}}
\put(120,5){\vector(1,0){50}}
\setcounter{c}{2}
\multiput(122,0)(5,0){2}{
   \makebox(0,0)[b]{\arabic{c}}
   \addtocounter{c}{1}
}
\put(138,0){\makebox{2+q}}
\multiput(123,5)(5,0){8}{\line(0,1){1}}
\put(123,18.5){\circle*{1}}
\put(128,18.5){\circle{1}}
\put(128,18.5){\line(3,1){15}}
\put(143,23.5){\circle*{1}}
\put(148,23.5){\circle{1}}
\end{picture}
\caption{$E_2^{(s,t)}(\HPbar)$ for $p$ odd.} \label{v1odd} 
\end{center}
\end{figure}
\begin{thm}[Adams, Mahowald, Miller, Harper-Miller]
\label{v1htp}
Consider Figures \ref{v1even},\ref{v1odd}.
\begin{enumerate}[(a)]
\item Above the lines and classes shown, the $E_2$ term vanishes.
\item The classes in dimensions $qk+3$ with filtration $k-1$ for $p$ odd and the 
classes in dimension $8k+7$ with filtration $\geq 4k+1$ for $p=2$ 
detect elements in $\pi_*\HPbar$ which under the isomorphism
$\pi_*\HPbar \to \pi_{*-1}S^3$ correspond to stable classes which are 
detected by the (stable, real) Adams $e$-invariant\footnote{See 
\eqref{defneinv} below for the definition of the $e$-invariant.}.
\end{enumerate}
\end{thm}
\begin{proof}
Recall from \cite{Ma1} and \cite{HM1} that for each prime there is a 
bigraded complex $(\Lambda(3),d)$ such that 
\[ E_2^{s,t}(S^3)= H^{s,t-3}\Lambda(3) \]
There is a short exact sequence of complexes
\[ 0 \lra \Lambda(1) \lra \Lambda(3) \lra W(1) \lra 0 \]
which induces a split short exact sequence on homology. $\Lambda(1)$ is a
complex with 0 differential which corresponds to the $\Z$-tower in $t-s=3$.

Let $M\otimes\Lambda$ denote the $E_1$-term of the stable Adams 
spectral sequence for the $\bmod  p$ Moore spectrum given by the 
$\Lambda$-algebra. The main results of \cite{Ma1} and \cite{HM1} state
 that there is a map of complexes $W(1) \lra M\otimes \Lambda$
 inducing an isomorphism on homology for 
$t-s<5s-16$ for $p=2$ and $t-s<(p+1)qs-(p+2)q$ for $p$ odd. This reduces
(a) to a statement about the stable $E_2$-term of the Moore spectrum
except in low dimensions where it is easily checked directly by computing
the homology of the complex $\Lambda(3)$.

The $E_2$-term of the Moore spectrum is computed in the required range in
\cite{Mi} for $p$ odd and in \cite{Ma3} for $p=2$. It agrees with the 
description given in Figures \ref{v1even} and \ref{v1odd}. In view of \eqref{relE2}, this 
completes the proof of (a).

In order to prove (b), let $M^{n}$ denote the $\bmod \ p$ Moore space with top cell in dimension
$n$ and $i:S^{n-1} \lra M^n$ denote the inclusion of the bottom cell.
Let $A:M^{n+r} \lra M^n$ with $r=8$ if $p=2$ and $r=q$ if $p>2$ be a map 
inducing an isomorphism in $K$-theory (see \cite[pp. 65-68]{Ad}).
Recall that, by construction, when $p=2$, the map 
$A$ is the composite of $4$ maps inducing $0$ on mod $2$
homology.

Suppose first that $p$ is odd. Let $\tilde{\alpha}$ denote an 
extension of a generator $\alpha \in \pi_{2p}S^3$ to $M^{2p+1}$. 
Then by \cite[Proposition 12.7]{Ad}, $\alpha_k=\tilde{\alpha}\circ A^{k-1} 
\circ i \in \pi_{qk+2}S^3$ survives to a stable class which is detected
by the $e$-invariant. The image of $\Sigma\alpha_k$ in 
$\pi_{qk+3} \HP^\infty$ under the inclusion $S^4 \hookrightarrow \HP^\infty$ 
lifts uniquely to a class $\overline{\Sigma \alpha_k} \in 
\pi_{qk+2} \HPbar$. Since $\alpha_k$ has Adams filtration at least $k$,
$\overline{\Sigma \alpha_k}$ has Adams filtration at least $k-1$ and so it must be represented
in $E_2(\HPbar)$ by the class in bidegree $(t-s,s)=(qk+3,k-1)$.

Next let $p=2$. 
Let $\mu \in \pi_{12}S^3$ be the unique element with Adams filtration $5$ (corresponding
to the leftmost class of the top lightning flash in Figure \ref{v1even}, for $k=1$). 
The stabilization map $\pi_{12} S^3 \to \pi_9S^0$ is injective (see \cite{To})
and can not decrease Adams filtration. It follows that $\mu$ suspends to the class 
$\mu_9 = \overline{\eta} \circ A \circ i$ in $\pi_9 S^0$ constructed in \cite[Theorem 12.13]{Ad}  
(there is only one class of filtration $5$ in the $9$-stem, see for instance \cite[p. 202]{MT}).
Let $\tilde \mu$ be an extension of $\mu$ to $M^{13}$ chosen to have Adams 
filtration $5$, which is possible as the map $S^3\{5\} \to S^3$ determined by a minimal
Adams resolution is injective on $\pi_{12}$ by Lemma \ref{coninj}(b). Then 
$\tilde\mu_k=\tilde{\mu}\circ A^{k-1}\circ i\in \pi_{8k+4}S^3$ has Adams filtration
 at least $4k+1$. Lemma \ref{coninj} implies that $S^3\{4k+1\} \to S^3$
 and $S^3\{4k+2\} \to S^3\{4k+1\}$ are injective on $\pi_{8k+4}$ from which it is easy to check that
 the Toda bracket $\langle \tilde \mu_k,2,\eta \rangle \in \pi_{8k+6}S^3$ has Adams filtration
at least $4k+2$. \cite[Proposition 12.18]{Ad} shows that this bracket
survives to a stable class of order $4$ detected by the $e$-invariant so
the bracket must be represented on $E_2(S^3)$ by the class in bidegree 
$(t-s,s)=(8k+6,4k+2)$. As in the case when $p$ is odd, this bracket corresponds
to an element in $\pi_{8k+7}\HPbar$ with Adams filtration $4k+1$ represented
in $E_2(\HPbar)$ by the class in bidegree $(t-s,s)=(8k+7,4k+1)$. This completes the proof of (b).
\end{proof}

\section{Construction of maps}
\label{sectionmaps}

Let $\nu: S^{4n+3} \to \HP^n$ denote the Hopf map. The cofiber sequence
\[ S^{4n+3} \llra{\nu} \HP^n \to \HP^{n+1} \]
shows that the obstruction to extending a map $f: \HP^n \to \HP^\infty$ 
to $\HP^{n+1}$ is 
\begin{equation}
\label{obst}
 o(f):= [f \circ \nu ] \in \pi_{4n+3} \HP^\infty.
\end{equation}
Alternatively, one can regard the obstruction as an element in $\pi_{4n+3} \HP^{n+1}$,
in which case, it is given by the formula 
\begin{equation}
\label{obst2}
o(f):= [f \circ \nu] - \deg(f)^{n+1} [\nu] \in \pi_{4n+3} \HP^{n+1} 
\end{equation}
as a simple cohomology calculation shows.
Note that, as the equivalence $S^3 \to \Omega \HPi$ factors through $\Omega S^4$, any element of $\pi_{*} \HPi$ 
(in particular the obstruction to extension) factors canonically through the bottom cell 
of $\HP^\infty$. We will still write $o(f)$ for the corresponding element in $\pi_{4n+3} S^4$.

We will write $KO(X), K(X), KSp(X)$ for the reduced $K$-groups of a space $X$. The
Atiyah-Hirzebruch spectral sequence implies that $K^0(\HP^n)=\Z[x]/x^{n+1}$.
In \cite{FG}, Feder and Gitler show that the Adams operations on $K(\HP^n)$ are determined by the formula
\begin{equation}
\label{adamsformula}
\psi^l(x) = l^2 x + \ldots + \frac{2l^2(l^2-1) \ldots (l^2-(n-1)^2)}{(2n)!} x^n.
\end{equation}
Moreover, they show that $\phi \colon K^0(\HP^n) \otimes \R \to K^0(\HP^n) \otimes \R$ is a map 
of rings commuting with Adams operations iff $\phi=\psi^l$ (i.e. $\phi$ is determined by \eqref{adamsformula})
for some real or purely imaginary number $l$. In particular, given a self map $f$ of $\HP^n$, the induced
homomorphism on $K$-theory is $f^*=\psi^{\sqrt{\deg(f)}}$.

Let $\cat A$ denote Adams' category of finitely generated abelian groups with real Adams 
operations  \cite[Section 6]{Ad}. 
\begin{lemma}
\label{extA}
Let $f$ be a self map of $\HP^n$ of degree $k$. Then $k \in \mathbf{FG}_{n+1}$ if and only if there is a 
self map $\phi$ of $KO(\Sigma^4\HP^{n+1})$ in $\A$ such that the following diagram commutes.
\begin{equation}
\label{parap}
\xymatrix{
  0 \ar[r] & KO(S^{4n+8}) \ar[r] \ar[d]^{k^{n+1}} & KO(\Sigma^4\HP^{n+1})
  \ar[r] \ar[d]^\phi  & KO(\Sigma^4\HP^n) \ar[r] \ar[d]^{\Sigma^4 f^*} & 0\\
  0 \ar[r] & KO(S^{4n+8}) \ar[r] & KO(\Sigma^4\HP^{n+1}) \ar[r] & 
 KO(\Sigma^4\HP^n) \ar[r] & 0 
}
\end{equation}
\end{lemma}
\begin{proof}
As $\psi^{-1}$ acts trivially, $K(\HP^n)$, is an object of $\A$. For each real or purely imaginary number
$l$, let $\psi^l_n \colon K(\HP^n)\otimes \R \to K(\HP^n)\otimes \R$ denote the map of real vector spaces
with Adams operations determined by \eqref{adamsformula} and multiplicative extension.
If $\varphi$ is a map of vector spaces such that 
\begin{equation}
\label{parapR}
\xymatrix{
  0 \ar[r] & K(S^{4n+4})\otimes \R \ar[r] \ar[d]^{k^{n+1}} & K(\HP^{n+1})\otimes \R
  \ar[r] \ar[d]^\varphi  & K(\HP^n) \otimes \R \ar[r] \ar[d]^{f^*=\psi^{\sqrt k}_n} & 0\\
  0 \ar[r] & K(S^{4n+4})\otimes \R \ar[r] & K(\HP^{n+1})\otimes \R \ar[r] & 
 K(\HP^n) \otimes \R\ar[r] & 0 
}
\end{equation}
commutes, then we can write (using the monomial basis for $K(\HP^{n+1})$) 
$$ \varphi =  \left[
\begin{matrix} \psi_n^{\sqrt k} & 0 \\ w & k^{n+1} \end{matrix}\right]  = \psi_{n+1}^{\sqrt k} +  \left[
\begin{matrix} 0 & 0 \\ v & 0 \end{matrix}  
\right]
$$
for some $v \in \R^n$. The map $\varphi$ commutes with $\psi^l_{n+1}$ iff $v \psi^l_n = l^{2n+2} v$.
But in view of \eqref{adamsformula}, $l^{2n+2}$ is not an eigenvalue of $\psi^l_n$ (except for $l \in \{-1,0,1\}$)
and hence $\varphi = \psi^{\sqrt k}_{n+1}$. In particular, $\varphi$ must be a ring homomorphism.
Conversely, $\varphi = \psi^{\sqrt k}_{n+1}$ gives a map making \eqref{parapR} commute.

The integer $k$ belongs to  $\mathbf{FG}_{n+1}$ iff
 $\varphi(x)$ is in the image of the restriction map $KSp(\HP^n) \to
K(\HP^n)$ (i.e. $\varphi(x)=a_1 x + \ldots a_n x^n$ with $a_i \in \Z$ and $a_{2i+1} \in 2\Z$).
Since  $\varphi$ is a ring homomorphism this is equivalent to saying that
the matrix $[a_{ij}]$ representing $\varphi$ in the monomial basis has 
\begin{equation}
\label{integcond}
a_{ij} \in \Z \quad \text{ and } \quad a_{(2i+1)(2j+1)} \in 2\Z.
\end{equation}

Given $M$ in $\A$, let $\Sigma^4 M$ denote the object in $\A$ with underlying abelian group $M$ and
 $\psi^l_{\Sigma^4M}=l^2\psi^l_M$. Note that when $M$ is torsion free, the self maps
of $M$ and $\Sigma^4M$ in $\A$ are the same. There is a natural inclusion in $\A$ 
$$ 
\xymatrix{ 
\Sigma^4K(\HP^n) \ar@{^{(}->}[r]^j & KO(\Sigma^4\HP^n)
} $$
determined by the zigzag
$$
\xymatrix{
\Sigma^4K(\HP^n) & \ar@{_{(}->}[l]_{\Sigma^4c} \Sigma^4KO(\HP^n) \ar@{^{(}->}[r]^p & KO(\Sigma^4 \HP^n)
}$$
where $c$ denotes the complexification map and $p$ is given by multiplication with the
generator of $KO^{-4}(*)$. The map $j\otimes \R$ is an isomorphism of vector spaces with Adams operations 
and hence, upon tensoring with $\R$,  the existence of a map $\phi\otimes \R$ making \eqref{parap}$\otimes\R$
commute is guaranteed. Moreover we must have 
$\phi \otimes \R = (j \otimes \R) \circ  \varphi \circ (j\otimes \R)^{-1} $.  

One easily checks that the integrality condition \eqref{integcond} on the matrix representation of $\varphi$ 
precisely corresponds to the integrality of $\phi \otimes \R$, i.e. to the
 existence of the map $\phi$ in the 
statement. This completes the proof.
\end{proof}
Recall from \cite{Ad} that the (real) $e$-invariant is a group homomorphism 
\begin{equation}
\label{defneinv}
 Z  \llra{e}  \Ext_{\cat A}(KO(X),KO(S^{j+1}))
\end{equation}
where $Z=\{\alpha \in \pi_j(X) | KO(\alpha)=0\}$. $e(\alpha)$ is defined to be the 
Yoneda class of the extension
\[ 0 \ra KO(S^{j+1}) \ra KO(X\cup_\alpha e^{j+1})\ra KO(X) \ra 0. \]
If $X$ is a sphere there is a natural identification of the target of $e$ with
 a subgroup of $\Q/\Z$. The stable $e$-invariant of $\beta \in \pi_{k+l}S^k$ is defined by suspending
$\beta$ to a class $\beta' \in \pi_{8m+l}S^{8m}$ for some $m$ and then setting
$e^s(\beta)=e(\beta') \in \Q/\Z$. This is independent of the
choice of $m$. 
\begin{prop}
\label{einv}
Let $f$ be a self map of $\HP^n$ and let $o(f)\in \pi_{4n+3} S^4$ denote its obstruction 
to extension. Then $\deg(f) \in \mathbf{FG}_{n+1}$ if 
and only if the stable $e$-invariant of $o(f)$ vanishes.
\end{prop}
\begin{proof}
Let $k=\deg(f)$. By Lemma \ref{extA}, $k \in \mathbf{FG}_{n+1}$ iff
the map of extensions \eqref{parap} exists. Writing $\circ$ for the Yoneda 
product in the category of extensions,  and $\Sigma^4\nu:S^{4n+3} \to \Sigma^4 \HP^n$ 
for the attaching map of the top cell of $\Sigma^4\HP^{n+1}$, 
the existence of \eqref{parap} is translated into the following equalities in $\Ext_{\cat A}$:
\begin{equation*}
\begin{split}
e(\Sigma^4\nu) \circ \Sigma^4 f^* & =   k^{n+1} \circ e(\Sigma^4\nu)  \\
e(\Sigma^4(f \circ \nu)) &  =  e(k^{n+1}\Sigma^4\nu) \\
e(\Sigma^4 o(f)) & =   0  \in \Ext_{\cat A}(KO(\Sigma^4HP^n), KO(S^{4n+8}))
\end{split}
\end{equation*}
where we have used that the $e$-invariant sends compositions of maps to Yoneda products, is a group homomorphism and \eqref{obst2}.
Since $o(f)$ factors through the inclusion of
 the bottom cell of $\HP^n$, letting $X=\Sigma^4\HP^n\cup_{o(f)}e^{4n+8}$, 
we have a map of extensions
\[
\xymatrix{
  0 \ar[r] & KO(S^{4n+8}) \ar[r] \ar[d]^= & KO(X)
  \ar[r] \ar[d]  & KO(\Sigma^4\HP^n) \ar[r] \ar[d]^{i^*} & 0\\
  0 \ar[r] & KO(S^{4n+8}) \ar[r] & KO(S^8\cup_{\Sigma^4o(f)}e^{4n+8})
 \ar[r] & KO(S^8) \ar[r] & 0 
}
\]
where $i$ denotes the inclusion of the bottom cell. We conclude that $e^s(o(f)) \circ i^* = e(\Sigma^4o(f))$.
The long exact sequence in $\Ext_{\cat A}$ induced by the short exact sequence 
\[ 0 \lra KO(\Sigma^4(\HP^n/\HP^1)) \lra KO(\Sigma^4\HP^n) \llra{i^*}
 KO(S^{8}) \lra 0 \]
together with the easily checked fact that $\Hom_\cat{A} 
(KO(\Sigma^4(\HP^n/\HP^1)), KO(S^{4n+8}))=0$ show that composition with
$i^*$ is injective which concludes the proof.
\end{proof}

We are now ready to prove our main result. Consider the function
\begin{equation}
\label{defnf}
f(p,n) = \begin{cases}
e(p,n) & \text{ for } (p=2, n\leq 5) \text{ or }  (p=3, n\leq 6) \\
&  \quad  \text{ or } \left(p>3, n \leq 1+ \frac{(p-1)^2}{2} \right), \\
4\lfloor \frac{n}{2}\rfloor -3 & \text{ for } p=2, n>5, \\
\lceil \frac{2n-2}{p-1}\rceil -1 & \text{ otherwise.} \\
\end{cases}
\end{equation}
where $\lceil x \rceil$ denotes the least integer $\geq x$. Theorem \ref{mainthm} is now an 
immediate consequence of the following result.
\begin{thm}
\label{main2}
$p^{f(p,n)}\Z_{(p)} \subset \mathbf{R}_{n,p}$.
\end{thm}

For a fixed $p$,  the function $e(p,n)$ (defined in \eqref{exponent}) grows logarithmically while $f(p,n)$ grows linearly.
Thus there remains a large gap between Theorem \ref{main2} and Conjecture \ref{conje}.
For the convenience of the reader, we include Table \ref{tab} comparing the values of $e(p,n)$ and $f(p,n)$ for 
small values of $p$ and $n$.

\begin{table}[h]
\label{tab}
\begin{tabular}{ c || c | c || c | c || c | c}
$n$ & $e(2,n)$ & $f(2,n)$ & $e(3,n)$ & $f(3,n)$ & $e(5,n)$ & $f(5,n)$ \\ \hline \hline
1 & 0 & 0 & 0 & 0 & 0 & 0  \\ \hline
2 & 3 & 3 & 1 & 1 & 0 & 0  \\ \hline
3 & 3 & 3 & 2 & 2 & 1 & 1 \\ \hline 
4 & 5 & 5 & 2 & 2 & 1 & 1 \\ \hline
5 & 5 & 5 & 2 & 2 & 2 & 2 \\ \hline
6 & 5 & 9 & 3 & 3 & 2 & 2 \\ \hline
7 & 5 & 9 & 3 & 5 & 2 &  2 \\ \hline
8 & 7 & 13 & 3 & 6 & 2 &  2 \\ \hline
9 & 7 & 13 & 4 & 7 & 2 &  2 \\ \hline
10 & 7 & 17 & 4 & 8 & 2 & 4  \\ \hline
\end{tabular}
\smallskip
\caption{Comparison of $e(p,n)$ and $f(p,n)$ for $p\leq 5$ and $n\leq 10$ .}\label{tab}
\end{table}

The two inputs needed for the proof of Theorem \ref{main2} are the (local versions of) Proposition \ref{einv} 
and the information about the unstable Adams spectral sequence contained in Theorem \ref{v1htp}.

Using exclusively Proposition \ref{einv}, we see that the $p$-local version of Conjecture \ref{conje} holds for
$n$ as long as the obstruction groups $\pi_{4j-2} S^3_{(p)}$ are detected by the $e$-invariant for $j\leq n$. 
As a consequence of the work of Toda \cite{To} and Adams \cite{Ad} (precise references will be given in the 
proof below) this happens in the following ranges: $n \leq 3$ for $p=2$, $n\leq 6$ for $p=3$ and 
$n\leq 1+ \frac{(p-1)^2}{2}$ for $p>3$. 
The range for $p>3$ can be roughly doubled using the techniques of this paper (cf. Remark \ref{finalrem}(a) below)
but we will not pursue this here. The strength of the $e$-invariant thus accounts
 for the first line in the definition of $f(p,n)$ except for the cases when $p=2$ and $n=4,5$.
 
The remaining cases in the definition of $f(p,n)$ come from putting together Proposition \ref{einv}
and  the description of the unstable Adams spectral sequence along the vanishing line in Theorem \ref{v1htp}. 
The slopes of the functions $n \mapsto f(p,n)$ are determined by the vanishing line of the unstable Adams 
spectral sequence but,  crucially for the cases when $p=2$ and $n=4,5$, the fact that the classes 
along the vanishing line are detected by the $e$-invariant effectively lowers the $y$-intercept of the vanishing line.
At odd primes, this ``lowering of the vanishing line'' is not enough to prove the $p$-local Feder-Gitler conjecture
 in any cases other than those settled by the method of the previous paragraph.
Nevertheless, it permits the construction of a number of maps for higher values of $n$.

\begin{proof}[Proof of Theorem \ref{main2}]
The obstruction to extending a map $\HP^{n-1} \to \HP^\infty_{(p)}$ lies in $\pi_{4n-1} \HPi_{(p)} = \pi_{4n-2} S^3_{(p)}$. 
For $p=2$, this group is detected by the $e$-invariant for $n\leq 3$. This is just the case when $p=2$ and $k=0$ of 
Theorem \ref{v1htp} (b). For $p=3$, the obstruction group is detected by the $e$-invariant for $n \leq 6$. This 
follows from the fact that $\pi_{4n-2} S^3_{(3)} = \Z/3$ for $2\leq n \leq 6$ (see 
\cite[Theorems 13.4,13.9,13.10]{To}) as Theorem \ref{v1htp} (b) gives a $\Z/3$ summand
 in these groups which is detected by the $e$-invariant.
Similarly, for $p>3$, the obstruction group is detected by the $e$-invariant for $n\leq 1+ \frac{(p-1)^2}{2}$ (but
see also Remark \ref{finalrem}(a)). Indeed, $\pi_{4n-2} S^3_{(p)} $ is either $\Z/p$ or $0$ in this range of values of $n$ by
\cite[Theorem 13.4]{To} and, when it is not $0$, Theorem \ref{v1htp} (b) provides a $\Z/p$ summand in the group
which is detected by the $e$-invariant.
By (the $p$-local version of) Proposition \ref{einv}, the proof of the Theorem 
when $p=2,n\leq 3$ or $p=3,n\leq 6$ or $p>3, n \leq 1+ \frac{(p-1)^2}{2}$ is complete.

First note that, in the remaining cases, $f(p,n)\geq e(p,n)\geq 2$.
Pick an Adams resolution for $\HPbar$ constructed from a minimal resolution. 
Then the Adams covers $\HPbar\{s\}$ are $3$-connected, 
$\pi_4(\HPbar\{s\})=\Z_{(p)}$ and the
maps $\HPbar\{s\} \to \HPi$ induce multiplication by $p^{s+1}$ on $\pi_4$. Thus, in order to prove the theorem, it 
is enough to provide a map $\HP^n \to \HPbar\{f(p,n)-1\}$ of degree $k$ on the bottom cell for 
an arbitrary $k \in \mathbb{Z}_{(p)}$. Such a map will be obtained by extending 
a map $\HP^1 \to \HPbar\{f(p,n)-1\}$ of degree $k$. 

Consider first the case when $p=2$ and suppose $n$ is even. Then for $k<n-1$ the obstruction to extension 
of the map $\HP^k \to \HPbar\{f(2,n)-1\}$ to $\HP^{k+1}$ lies in the group $\pi_{4k+3} \HP^\infty\{f(2,n)-1\}$.
By Lemma \ref{coninj}(a) and the vanishing line of Theorem \ref{v1htp}(a) this group
 is 0 so the map extends. Thus we have a map $\HP^{n-1} \to \HPbar\{f(2,n)-1\}$.
  Since the Hopf map $S^{4n-1} \to \HP^{n-1}$ has Adams
filtration $1$, the composite $S^{4n-1} \to \HPbar\{f(2,n)-1\}$ factors through $\HPbar\{f(2,n)\}$. Let
$\gamma:S^{4n-1} \to \HPbar\{f(2,n)\}$ denote this obstruction class. By Theorem \ref{v1htp}(b),
the image $\gamma'$ of $\gamma$ in $\pi_{4n-1} \HP^\infty$ is detected by the $e$-invariant (it is in the 
summand of order $4$ corresponding to one of the top lightning flashes). Since $\gamma'$ is
the obstruction to extension of a map $\HP^{n-1} \to \HP^\infty$ with degree in $\mathbf{FG}_{n,2}$ (because $f(p,n)
\geq e(p,n)$), Proposition \ref{einv} implies that $\gamma'=0$.
Since the map $\pi_{4n-1} \HPbar\{f(2,n)\} \to \pi_{4n-1} \HPi$ is injective by Lemma \ref{coninj}
(b), we conclude that $\gamma=0$ and hence our map extends to a map $\HP^n \to \HPbar\{f(2,n)-1\}$
of degree $k$. This completes the proof of the case when $p=2$ and $n$ is even. 

The case when $p=2$ and $n$ is odd follows, since $f(2,2m+1)=f(2,2m)$ and, for $n=2m+1$, the map
$\HP^{2m} \to \HPbar\{f(2,2m)-1\}$ described in the previous paragraph extends one more stage. Indeed,  
its obstruction to extension factors through $\HPbar\{f(2,2m)\}$ and hence vanishes by Lemma \ref{coninj}(a) 
(using the vanishing line of Theorem \ref{v1htp}(a)). This completes the proof for $p=2$.

For $p$ odd, the reasoning is identical: suppose $j\geq 2$ and we want to extend a map 
$\HP^1 \to \HPbar\{j-1\}$ of degree $k \in \Z_{(p)}$.
By Lemma \ref{coninj}(a) and the vanishing line of Theorem \ref{v1htp}(a) there is no obstruction
to extending it to $\HP^n$ for $n\leq \frac{j(p-1)}{2}$.
For $n= \frac{j(p-1)}{2}+1$ there is an obstruction but, since the Hopf map 
has Adams filtration $1$, the obstruction factors through $\HPbar\{j\}$ and hence 
vanishes. By Theorem $\ref{v1htp}$ there are no further obstructions until $n= \frac{(j+1)(p-1)}{2}+1$.
At that point, there is an obstruction, which again factors through $\HPbar\{j\}$. This obstruction 
is detected by the $e$-invariant and hence, by Proposition \ref{einv}, vanishes as long as 
$j\geq e\left(p, \frac{(j+1)(p-1)}{2}+1\right)$. This last condition holds because $j\geq 2$ and
hence we obtain a map $\HP^{\frac{(j+1)(p-1)}{2}+1} \to \HPi_{(p)}$ of degree $kp^j$.

For a given $n$, the least $j$ such that $n \leq \frac{(j+1)(p-1)}{2} + 1$ is $\lceil \frac{2n-2}{p-1}\rceil -1$.
In the range of values of $n$ under consideration (namely $n\geq 7$ for $p=3$ and 
$n> 1+ \frac{(p-1)^2}{2}$ for $p>3$) we have $f(p,n)= \lceil \frac{2n-2}{p-1}\rceil -1 \geq 2$. 
We conclude that there are maps of degree $kp^{f(p,n)} \colon \HP^n \to \HPi_{(p)}$ for all $k \in \Z_{(p)}$.
This completes the proof.
\end{proof}

\begin{remark}
\label{finalrem}
For $p>3$, Theorem \ref{main2} can be improved using the following observations:
\begin{enumerate}[(a)]
\item For $n < \frac{(2p+1)(p-1)}{2}$ all classes in $\pi_{4n-2}(S^3_{(p)})$ are detected by 
the $e$-invariant, therefore the conjecture holds at odd primes for these values of $n$.
This follows from the results of Mahowald and Harper-Miller (cf. proof of Theorem \ref{v1htp})
 together with calculations with the $\Lambda$-algebra in
low degrees. See \cite[Corollary 2.9]{G} for a proof.
\item Let $g(p,n)= \sup \{ k : p^k \Z_{(p)} \subset \mathbf{R}_{n,p} \}$. 
Using \eqref{obst2}, it is easy to show that, for $f,g$ self maps of $\HP^n$, the obstruction to extension
of the composite is $o(f g) = \deg(g)^n o(f) + \deg(f) o(g)$ (see \cite[Proposition 4.6(b)]{G}). Since
(by Proposition \ref{rector}) there is a self map of $\HP^n_{(p)}$ of degree $p^2$, Selick's exponent theorem \cite{Se} 
implies that  
\begin{equation}
\label{bound}
g(p,n+1)\leq 2+g(p,n).
\end{equation}
For $p>3$, the first positive value of the function $n \mapsto f(p,n)-e(p,n)$ is
$p-3$ (with the improved range of Remark \ref{finalrem} (a), it would be $2p-3$). When this number is greater 
than $2$, we can apply \eqref{bound} to improve on Theorem \ref{main2} in a very limited range of values of $n$.
\end{enumerate}
\end{remark}

\begin{remark}
\label{finalfinalrem}
The first open case of Conjecture \ref{conje} is $n=6$. In view of Theorem \ref{main2}, the only prime
that needs to be considered is $p=2$. The relevant obstruction group is $\pi_{22} S^3_{(2)} = \Z/2 \oplus \Z/4$ (cf. \cite
[Theorem 12.9]{To}). As the $\Z/4$ summand corresponds to a class detected by the $e$-invariant, in order to 
settle the conjecture for $n=6$ it suffices to find an effective way of detecting the remaining $\Z/2$ summand.
This $\Z/2$ is generated by an unstable $v_1$-periodic class appearing in the bottom 
lightning flash in Figure \ref{v1even}. At an odd prime, the first value of $n$ for which the conjecture is open 
is $n=\frac{(2p+1)(p-1)}{2}$. The obstruction group for this value of $n$ is $\pi_{(2p+1)q-1}S^3_{(p)} = \Z/p \oplus \Z/p$.
One of the $\Z/p$ summands is detected by the $e$-invariant whilst the other is not $v_1$-periodic (see 
\cite[Corollary 2.9]{G} for more details).
\end{remark}


\begin{thebibliography}{MG2}
\bibitem[Ad]{Ad}
J. F. Adams, \emph{On the groups J(X).IV}, Topology 5 (1966) 21-71.
\bibitem[AC]{AC}
M. Arkowitz and C. R. Curjel, \emph{On maps of $H$-spaces}, Topology 6
 (1967) 137-148. 
\bibitem[CM]{CM}
E. Curtis and M. Mahowald,
 \emph{The unstable Adams spectral sequence for $S^3$}, Contemp. Math. 96,
AMS (1989) 125-162.
\bibitem[FG]{FG}
S. Feder and S. Gitler, \emph{Mappings of quaternionic projective spaces},
Bol. Soc. Mat. Mex. 18 (1973) 33-37.
\bibitem[GS]{GS} D. L. Gonçalves and M. Spreafico, \emph{Quaternionic line bundles over quaternionic projective spaces.} Math. J. Okayama Univ. 48 (2006), 87--101.
\bibitem[G]{G} G. Granja, \emph{On quaternionic line bundles}, Ph.D. Thesis, M.I.T., 1999.
\bibitem[Ha]{Ha} J. Harper, \emph{Secondary cohomology operations.} Graduate Studies in Mathematics, 49. American Mathematical Society, Providence, RI, 2002. xii+268 pp.
\bibitem[HM1]{HM1}
J. Harper and H. R. Miller, \emph{On the double suspension homomorphism
at odd primes}, Trans. Amer. Math. Soc. 273 (1982) 319-331.
\bibitem[HM2]{HM2}
J. Harper and H. R. Miller, \emph{Looping Massey-Peterson towers} in:
S.M. Salamon, B. Steer, W.A. Sutherland, eds. \emph{Advances in Homotopy
Theory}(London Math. Soc. Lec. Notes 139, Cambridge University Press,
Cambridge, 1989) 69-86.
\bibitem[IMO]{IMO} N. Iwase, K. Maruyama and S. Oka, \emph{A note on 
$\mathcal E(\HP^n)$ for $n\leq 4$.} Math. J. Okayama Univ. 33 (1991), 163--176. 
\bibitem[Ma1]{Ma3}
M. Mahowald, \emph{The order of the image of the $J$ homomorphism}, Bull. Amer. Math. Soc.
{\bf 76} (1970) 1310-1313.
\bibitem[Ma2]{Ma1}
M. Mahowald, \emph{On the double suspension homomorphism}, Trans. Amer.
Math. Soc. 214 (1975) 169-178.
%\bibitem[Ma2]{Ma2}
%M. Mahowald, \emph{The image of $J$ in the $EHP$ sequence}, Annals of Math.
%116 (1982) 65-112.
\bibitem[MP]{MP}
W. Massey and F. Peterson, \emph{The mod 2 cohomology structure of certain fibre spaces.}, Mem. Amer.
Math. Soc. {\bf 74} (1967).
\bibitem[MR]{MR}
H. Marcum and D. Randall, \emph{A note on self-mappings of quaternionic 
projective spaces}, An. Acad. Brasil. Ci. 48 (1976) 7-9.
\bibitem[MG1]{MG1}
C. McGibbon, \emph{Self maps of projective spaces}, Trans. Amer. Math. Soc.
 271 (1982) 325-346.
\bibitem[MG2]{MG2}
C. McGibbon, \emph{Multiplicative properties of power maps. II}, 
Trans. Amer. Math. Soc. 274 (1982) 479-508.
\bibitem[Mi]{Mi}
H. Miller, \emph{A localization theorem in homological algebra}, Math. Proc.
Camb. Phil. Soc. 84 (1978) 73-84.
\bibitem[Ms]{Ms}
G. Mislin, \emph{The homotopy classification of self-maps of infinite quaternionic projective space},
Quart. J. Math. Oxford (2) 38 (1987) 245-257.
\bibitem[MT]{MT}
R. Mosher and M. Tangora, \emph{Cohomology Operations and Applications in Homotopy Theory}, 
Harper and Row, Publishers, New York - London, 1968, x+ 214 pp. 
\bibitem[Re]{Re}
D. Rector, \emph{Loop structures on the homotopy type of $S^3$}, Lecture
Notes in Math. 249 (1971) 99-105. 
\bibitem[Se]{Se}
P. Selick, \emph{Odd primary torsion in $\pi_kS^3$}, Topology 17 (1978)
407-412.
\bibitem[St]{St}
J. Stasheff, \emph{$H$-spaces from a homotopy point of view}, Lect. Notes in Math. {\bf 161},
Springer, Berlin (1970).
\bibitem[Su]{Su}
D. Sullivan, \emph{Geometric Topology I. Localization, periodicity and Galois
symmetry}, MIT Notes (1970).
\bibitem[Ta]{Ta}
M. Tangora, \emph{Computing the homology of the lambda algebra}, Memoirs of
the Amer. Math. Soc. 337 (1985).
\bibitem[To]{To}
H. Toda, \emph{Composition methods in homotopy groups of spheres}, Annals of 
Mathematics Studies 49, Princeton University Press (1962). 
\end{thebibliography}
\end{document}